\documentclass[11pt,reqno]{amsart}
\usepackage[a4paper,margin=28mm]{geometry}
\usepackage{amsmath,amssymb,amsthm,mathtools,mathrsfs}
\usepackage{hyperref}
\usepackage{cleveref}

\newtheorem{theorem}{Theorem}

\newtheorem{lemma}[theorem]{Lemma}
\theoremstyle{remark}
\newtheorem{remark}[theorem]{Remark}

\newcommand{\g}{\mathfrak g}
\newcommand{\h}{\mathfrak h}
\newcommand{\C}{\mathbb C}

\newcommand{\Sym}{\operatorname{Sym}}
\newcommand{\grad}{\nabla}
\newcommand{\on}{\operatorname}

\newcommand{\W}{\mathscr{W}}

\title{Polynomial twisted periods and affine vertex algebras}
\date{}
\author{Tomoyuki Arakawa}
\address[Tomoyuki Arakawa]
{Okinawa Institute of Science and Technology\\
1919-1 Tancha, Onna-son, Kunigami-gun, Okinawa 904-0495, Japan}
\email{tomoyuki.arakawa@oist.jp}

\author{Wenbin Yan}
\address[Wenbin Yan]{
Yau Mathematical Sciences Center, Tsinghua University\\
C544 Shuangqing Complex Building, Haidian District, Beijing 100084, China
}
\email{wbyan@mail.tsinghua.edu.cn}

\begin{document}

\begin{abstract}
We prove the conjecture of Kang, Lawrie, and Song that Zhu's
$C_2$-algebra of a simply-laced boundary admissible affine vertex
algebra is isomorphic to the Jacobian ring of an invariant potential
obtained by lifting a polynomial twisted period from a Cartan
subalgebra via the Chevalley isomorphism. We further identify the
$C_2$-algebras of the corresponding Drinfeld--Sokolov reductions with
the Jacobian rings of the restrictions of this potential to Slodowy
slices.
We also extend the Jacobian presentation to non-simply-laced principal
boundary admissible affine vertex algebras, where the restriction of the
potential to the Cartan subalgebra satisfies a Dunkl-singularity equation with unequal
parameters.
\end{abstract}

\maketitle

\section{Introduction}

Polynomial twisted periods arise from the interplay between invariant
theory and the geometry of Frobenius manifolds. Saito, Yano, and
Sekiguchi \cite{SYS} introduced distinguished flat coordinates on the
orbit spaces of finite Coxeter groups, and Dubrovin \cite{Dub1} later
showed that these orbit spaces carry natural polynomial Frobenius
manifold structures. In his theory of almost duality, Dubrovin
\cite{Dub2} introduced twisted periods as deformed flat coordinates
for the dual structure. At special values of the deformation
parameter, these periods become polynomial. Feigin and Silantyev
\cite{FS} subsequently related such polynomial twisted periods to
singular polynomials in polynomial representations of rational
Cherednik algebras, thereby giving a representation-theoretic
realization of these objects in terms of Saito flat coordinates.

More recently, polynomial twisted periods have appeared in a rather
different context, namely, in the study of protected sectors of
four-dimensional $\mathcal N=2$ superconformal field theories. For the
Argyres--Douglas theories $\mathcal D_p(G)$, where $G$ is of ADE type
and $\gcd(p,h^\vee)=1$, Kang, Lawrie, and Song \cite{KLS} proposed that
the polynomial twisted period of degree $p+1$ and parameter
$p/h^\vee$, regarded via the Chevalley isomorphism as a $G$-invariant
polynomial $P$ on $\mathfrak g$, provides a universal potential for the
theory. The associated Jacobian ring
\[
R_P=\mathbb C[\mathfrak g]/(\partial P)
\]
is conjectured to encode OPE decoupling relations in the protected
sector. In this proposal, its reduced critical locus reproduces the Higgs branch, while
its nonreduced structure retains additional information not visible in
the Higgs-branch variety itself.

In view of the 4d/2d correspondence of Beem, Lemos, Liendo, Peelaers,
Rastelli, and van Rees \cite{BLLPRvR}, it is natural to expect this
additional information to be encoded in the corresponding vertex
algebra. A central prediction of this correspondence is that the
associated variety of the vertex algebra, introduced in \cite{Ara},
coincides with the Higgs branch of the corresponding 4d theory
\cite{BR}. By definition, the associated variety $X_V$ is the
underlying reduced scheme of the associated scheme
\[
\widetilde{X}_V:=\operatorname{Spec}R_V,
\]
where $R_V$ is Zhu's $C_2$-algebra \cite{Zhu}. Since $R_V$ need not
be reduced, it is natural to expect the additional nonreduced
structure of the protected sector to be captured by
$\widetilde{X}_V$. Kang, Lawrie, and Song \cite{KLS} conjectured that
this is indeed the case.

The purpose of the present paper is to prove this conjecture for all
simply-laced boundary admissible levels and to establish the
corresponding Jacobian presentation for Drinfeld--Sokolov reductions.

Let us state our main results more precisely. Let $\g$ be a
finite-dimensional simply-laced complex simple Lie algebra, let $G$
be its adjoint group, and let $h^\vee$ be its dual Coxeter number.
Fix a Cartan subalgebra $\h$ with Weyl group $W$. The chiral algebra
corresponding to $\mathcal D_p(G)$ was identified in \cite{SXY} with
the simple affine vertex algebra $L_k(\mathfrak g)$ at the boundary
admissible level \cite{KW1,KW2} determined by
\begin{equation}\label{eq:level}
k+h^\vee=\frac{h^\vee}{p}.
\end{equation}

\begin{theorem}[Kang--Lawrie--Song conjecture \cite{KLS}]
\label{th:main1}
Let $\mathfrak g$ be a simply-laced simple Lie algebra, let $p\geq2$
satisfy $\gcd(p,h^\vee)=1$, and let $k$ be given by
\eqref{eq:level}. Then Zhu's $C_2$-algebra
$R_{L_k(\mathfrak g)}$ is isomorphic, as a graded Poisson algebra, to
the Jacobian ring
\[
R_P=\mathbb C[\mathfrak g]/(\partial P),
\]
where $P\in\C[\g]^G$ is the Chevalley lift of a nonzero $W$-invariant
polynomial twisted period on $\h$ of degree $p+1$ and parameter
$p/h^\vee$. The Poisson bracket on $R_P$ is induced by the
Kirillov--Kostant bracket on $\C[\g]$. Under the assumption
$\gcd(p,h^\vee)=1$, such a polynomial twisted period is unique up to
multiplication by a nonzero scalar \cite{FS}.
\end{theorem}

In other words, the defining relations of the $C_2$-algebra of a
boundary admissible affine vertex algebra arise as the components of the gradient of a single
invariant polynomial, while the twisted-period equation itself emerges
from the affine singular vector through the radial identity
\eqref{eq:radial} below.

By \Cref{th:main1}, the associated variety $X_{L_k(\g)}$ coincides
with the reduced critical locus $\on{Crit}(P)$ of $P$. Hence, by \cite{Ara},
\[
\on{Crit}(P)=\{x\in\g\mid \on{ad}(x)^{2p}=0\}.
\]
In types $A$ and $D$, this description was obtained directly in
\cite[Theorem~1.2]{Kang}.

We next consider the quantized Drinfeld--Sokolov reduction
\cite{FF,KRW}
\[
H_{DS,f}^0(L_k(\mathfrak g))
\]
of $L_k(\mathfrak g)$ associated with a nilpotent element
$f\in\mathfrak g$.
The geometry of nilpotent Slodowy slices plays a central role in the
study of affine $W$-algebras at admissible levels; see, for example,
\cite{AvEM}. The compatibility of Drinfeld--Sokolov reduction with
the $C_2$-quotient \cite[Corollary~4.5.10 and Theorem~4.6.1]{Ara}
gives an isomorphism of Poisson schemes
\begin{equation}\label{eq:DS-scheme}
\widetilde{X}_{H_{DS,f}^0(L_k(\g))}
\cong
\widetilde{X}_{L_k(\g)}\times_{\g^*}\mathcal S_f,
\end{equation}
where $\mathcal S_f$ is the Slodowy slice through $f$, viewed as a
subvariety of $\g^*$ via the invariant form. The Poisson structure on
the right-hand side is the transverse Poisson structure induced by
Hamiltonian reduction; see, for example, \cite{GG}. Together with the
transversality argument in \Cref{lem:slice}, this yields the following
result.

\begin{theorem}\label{th:main2}
Let $\g$, $k$, and $P$ be as in \Cref{th:main1}, and let
$f\in\g$ be a nilpotent element in $\operatorname{Crit}(P)$.
Set $P_f=P|_{\mathcal S_f}$. Then
\[
R_{H_{DS,f}^0(L_k(\g))}
\cong
\C[\mathcal S_f]/(\partial P_f)
\]
as Poisson algebras, with the bracket on the right induced by the
Slodowy Poisson structure. Here $(\partial P_f)$ denotes the ideal
generated by the directional derivatives of $P_f$ along
$\mathcal S_f$.
\end{theorem}

The hypothesis on $f$ ensures that the reduction is nonzero
\cite[Proposition~4.6.2]{Ara}. It is conjectured to be simple, or,
equivalently, to be isomorphic to the simple quotient $\W_k(\g,f)$
of the universal $W$-algebra $\W^k(\g,f)$. This is known when
$f$ admits a good even grading: see \cite[Theorem~7.8]{AvE} when
$f$ lies in the dense orbit of $X_{L_k(\g)}$, and
\cite[Theorem~10.1]{ACF} for arbitrary $f\in X_{L_k(\g)}$.
Consequently, in this case \Cref{th:main2} gives
\[
R_{\W_k(\g,f)}
\cong
\C[\mathcal S_f]/(\partial P_f).
\]
In particular, this applies when $f$ is principal, in which case
simplicity was established in \cite{A}. The Jacobian presentation in
the principal case was conjectured earlier by Xie and Yan \cite{XY}.

Under the expected 4d/2d correspondence, the reduction
$H_{DS,f}^0(L_k(\mathfrak g))$ is associated with the nilpotently
Higgsed theory $\mathcal D_p(G,G\mathord{.}f)$ obtained from
$\mathcal D_p(G)$. From the viewpoint of polynomial twisted periods,
this Higgsing corresponds to restricting the potential $P$ to the
Slodowy slice $\mathcal S_f$.

When $f$ is principal and $p>h^\vee$, \Cref{th:main2} has a
particularly interesting interpretation from the viewpoint of the
4d/2d correspondence. The restriction $P_f$ of the invariant
potential to the principal Slodowy slice is closely related to the
Seiberg--Witten geometry of the corresponding Argyres--Douglas theory
\cite{KLS}. Our theorem therefore identifies the Jacobian algebra
naturally associated with this four-dimensional geometry with Zhu's
$C_2$-algebra of the corresponding simple principal $W$-algebra,
which is rational by \cite{Arational}. This provides an algebraic link
between the proposed Seiberg--Witten description and the associated
scheme of the chiral algebra and is compatible with the perspective
on four-dimensional mirror symmetry proposed by Shan, Xie, and Yan
\cite{SXY2}.

Finally, our argument extends to non-simply-laced principal boundary
admissible affine vertex algebras. In this case, we again obtain a Jacobian
presentation of Zhu's $C_2$-algebra, while the restriction of the potential to
the Cartan subalgebra satisfies a Dunkl-singularity condition with unequal
parameters; see \Cref{th:non}.

\medskip
\noindent
\textbf{Acknowledgements.}
Wenbin Yan is supported by the National Key R\&D Program of China
(Grant No.~2025YFA1017400).
Tomoyuki Arakawa is partially supported by JSPS KAKENHI Grant Numbers
25K21659 and 26H01997.
This paper was written during his stay at SLMath for the program
``Representation Theory Under the Influence of Quantum Field Theory.''
He thanks the institute for its excellent working environment.

\medskip
\noindent
\textbf{Declaration of AI use.}
During the preparation of this work, the authors used OpenAI's ChatGPT
and Microsoft 365 Copilot to assist with mathematical discussions,
the development and exposition of proofs, literature searches, and
the drafting and editing of the manuscript. The authors independently
verified all mathematical arguments and references, revised all
AI-assisted material, and take full responsibility for the content of
the paper.

\section{Polynomial twisted periods}\label{sec:periods}

Let $\g$ be a finite-dimensional simply-laced complex simple Lie
algebra, let $G$ be its adjoint group, and let $h^\vee$ be its dual
Coxeter number. Normalize the invariant bilinear form
$(\cdot,\cdot)$ so that every root has squared length two. We use this
form to identify $\g$ with $\g^*$, and similarly identify a Cartan
subalgebra $\h$ with $\h^*$. For a root $\alpha$, its image under this identification is the
coroot $\alpha^\vee$. Choose
root vectors satisfying
\[
(e_\alpha,f_\alpha)=1,
\qquad
[e_\alpha,f_\alpha]=\alpha^\vee.
\]

Fix an integer $p\geq2$ with $\gcd(p,h^\vee)=1$, and set
\[
k=-h^\vee+\frac{h^\vee}{p},
\qquad
\kappa=k+h^\vee=\frac{h^\vee}{p}.
\]

Choose a set of positive roots, and write
$\h_{\mathrm{reg}}$ for the complement of the root hyperplanes. A
nonzero homogeneous $W$-invariant polynomial
$q\in\C[\h]^W$ of degree $p+1$ is called a polynomial twisted period
with parameter $p/h^\vee$ if
\begin{equation}
\label{eq:twisted}
H_q(x)u
=
\frac{p}{h^\vee}
\sum_{\alpha>0}
\frac{\alpha(u)\alpha(\grad q(x))}{\alpha(x)}
\alpha^\vee,
\qquad
x\in\h_{\mathrm{reg}},\quad u\in\h.
\end{equation}
Here the Hessian and gradient are taken with respect to the fixed
bilinear form.

Equation~\eqref{eq:twisted} is the coordinate-free form of
\cite[Equation~(2.3)]{FS}, with structure constants given by
\cite[Equation~(2.1)]{FS} and parameter $\nu=p/h^\vee$.

For the parameters under consideration, such a polynomial exists and
is unique up to multiplication by a nonzero scalar by the
classification of Feigin--Silantyev \cite[Theorem~4.8]{FS}. For
clarity, we reserve the term \emph{polynomial twisted period} for $q$
on $\h$ and call its lift $P$ the \emph{invariant potential}. Thus the
potential $P$ in \Cref{th:main1} is the Chevalley lift of $q$, namely,
the element corresponding to $q$ under the Chevalley isomorphism
\[
\C[\g]^G\cong\C[\h]^W.
\]

For a polynomial $Q$ on an affine space, write $J_Q=(\partial Q)$ for
the ideal generated by its directional derivatives.

\section{Proof of the affine presentation}\label{sec:affine}

Write $I_k$ for the maximal proper graded ideal of the universal
affine vertex algebra $V^k(\g)$ associated with $\g$ at level $k$, so
that
\[
L_k(\g)=V^k(\g)/I_k.
\]
Let
\[
\sigma:V^k(\g)\longrightarrow R_{V^k(\g)}
=\Sym\g=\C[\g]
\]
be the $C_2$-quotient map, where
$\sigma(a_{(-1)}\mathbf1)(x)=(a,x)$. Recall that
\[
R_V=V/C_2(V),
\qquad
C_2(V)=\operatorname{span}\{a_{(-2)}b:a,b\in V\}.
\]

We write $a(n)=a_{(n)}$ for the affine current modes. The key
ingredient in the proof is the radial identity of
\cite[Lemma~2.1]{ADS}:
\begin{equation}
\label{eq:radial}
\kappa\,\partial_u\sigma(w)(x)
=
\sum_{\alpha>0}
\frac{\alpha(u)}{\alpha(x)}
\sigma\bigl(f_\alpha(0)e_\alpha(0)w\bigr)(x)
+
\sigma\bigl(u(1)w\bigr)(x),
\end{equation}
for $w\in V^k(\g)$, $u\in\h$, and
$x\in\h_{\mathrm{reg}}$.

\begin{proof}[Proof of \Cref{th:main1}]
\noindent
\emph{Step 1: the generating space is an adjoint representation.}
Let $\theta$ be the highest root and put $\lambda=k\Lambda_0$. The
simple roots of the integral root system of $\lambda$ may be chosen
to be the finite simple roots together with $\beta=p\delta-\theta$.
In the simply-laced normalization, $\beta^\vee=pK-\theta^\vee$, and
hence
\[
\langle\lambda+\widehat\rho,\beta^\vee\rangle
=
p(k+h^\vee)-\langle\rho,\theta^\vee\rangle
=
h^\vee-(h^\vee-1)
=
1.
\]
The singular-vector generation theorem for admissible highest-weight
modules \cite{KW0,KW1}, applied to the vacuum quotient, implies that
the maximal ideal $I_k$ is generated, as an affine Lie algebra module,
by a singular vector $v$ of weight
\[
s_\beta\mathbin{\circ}\lambda
=
\lambda-\beta
=
\lambda+\theta-p\delta.
\]
The submodules corresponding to the finite simple roots have already
been factored out in passing from the Verma module to the vacuum
module. Thus $v$ has conformal degree $p$ and finite highest weight
$\theta$.

The finite-dimensional $\g$-module
\[
E=U(\g(0))v
\]
is therefore isomorphic to the adjoint representation. Moreover,
\[
(\g\otimes z\C[z])E=0,
\qquad
I_k=U(\g\otimes z^{-1}\C[z^{-1}])E.
\]
Choose a $\g$-module isomorphism
\[
\g\longrightarrow E,\qquad a\longmapsto v_a.
\]
Define a homogeneous polynomial map $F:\g\to\g$ of degree $p$ by
\begin{equation}
\label{eq:F}
(F(x),a)=\sigma(v_a)(x).
\end{equation}
The equivariance of $\sigma$ implies that $F$ is $G$-equivariant.

\medskip
\noindent
\emph{Step 2: the components of $F$ generate the entire $C_2$-ideal.}
The preceding PBW description of $I_k$ gives
\begin{equation}
\label{eq:ideal}
\sigma(I_k)
=
\C[\g]\,\sigma(E)
=
\bigl((F(-),a)\mid a\in\g\bigr).
\end{equation}
Indeed, a mode $a_{(-r)}$ with $r\geq2$ has zero image in the
$C_2$-quotient, whereas $a_{(-1)}$ acts there by multiplication by
the linear function $(a,-)$. Reordering the negative modes produces
commutator terms whose mode numbers are at most $-2$, which therefore
vanish in the $C_2$-quotient. Conversely, every polynomial multiple
of an element of $\sigma(E)$ is the image of an iterated $(-1)$-mode
action on $E$.

We also have $F\neq0$: since $k\neq-h^\vee$, the image of the nonzero
singular vector $v$ in the $C_2$-quotient is nonzero by
\cite[Theorem~1.2]{AJM}. Thus this step does not require the known
classification of the associated variety.

\medskip
\noindent
\emph{Step 3: the radial identity \eqref{eq:radial} becomes the
twisted-period equation.}
Equivariance implies that $F(x)\in\h$ for $x\in\h_{\mathrm{reg}}$;
by polynomiality, the same holds for all $x\in\h$. Put $f=F|_\h$.
For $u,a\in\h$, apply \eqref{eq:radial} to $w=v_a$. The positive-mode
term vanishes because $u(1)E=0$, and
\[
f_\alpha(0)e_\alpha(0)v_a
=
v_{[f_\alpha,[e_\alpha,a]]}
=
\alpha(a)v_{\alpha^\vee}.
\]
We obtain
\[
\kappa(\partial_uf(x),a)
=
\sum_{\alpha>0}
\frac{\alpha(u)\alpha(a)}{\alpha(x)}
\alpha(f(x)),
\]
or equivalently,
\begin{equation}
\label{eq:f-equation}
\partial_uf(x)
=
\frac{p}{h^\vee}
\sum_{\alpha>0}
\frac{\alpha(u)\alpha(f(x))}{\alpha(x)}
\alpha^\vee.
\end{equation}
Pairing the right-hand side with $a$ gives an expression symmetric
in $u$ and $a$. Hence the polynomial one-form $(f(x),dx)$ is closed:
\begin{equation}
\label{eq:closed}
(\partial_uf(x),a)
=
(\partial_af(x),u).
\end{equation}
Although derived on $\h_{\mathrm{reg}}$, this polynomial identity
holds on all of $\h$.

Set
\[
q_0(x)=\frac{1}{p+1}(x,f(x)).
\]
Since $f$ is homogeneous of degree $p$, \eqref{eq:closed} and Euler's
identity give
\[
\partial_u q_0(x)
=
\frac{(u,f(x))+(x,\partial_uf(x))}{p+1}
=
\frac{(u,f(x))+(\partial_xf(x),u)}{p+1}
=
(f(x),u).
\]
Thus $\grad q_0=f$. The polynomial $q_0$ is $W$-invariant and
homogeneous of degree $p+1$. It is nonzero: if $f$ vanished, then
equivariance and the density of the regular semisimple elements would
force $F=0$. Substituting $f=\grad q_0$ into
\eqref{eq:f-equation} shows that $q_0$ satisfies
\eqref{eq:twisted}.

By the uniqueness of the polynomial twisted period \cite{FS},
$q_0=cq$ for some $c\in\C^\times$. In the present situation, the
possible repeated-exponent exception in even type $D$ does not occur:
the repeated exponent is $n-1$ for $D_n$, whereas
$h^\vee=2n-2$, so its residue class modulo $h^\vee$ is not coprime to
$h^\vee$.

\medskip
\noindent
\emph{Step 4: identify the ideal on the whole Lie algebra.}
Let $P_0=cP$ be the Chevalley lift of $q_0$. For regular $x\in\h$,
invariance implies that $\grad P_0(x)\in\h$, while restriction of the
differential gives
\[
\grad P_0(x)
=
\grad q_0(x)
=
f(x)
=
F(x).
\]
Both $\grad P_0$ and $F$ are $G$-equivariant polynomial maps. Their
equality on $\h_{\mathrm{reg}}$ therefore implies their equality on
the dense set of regular semisimple elements, and hence on all of
$\g$. Equation~\eqref{eq:ideal} now yields
\[
\sigma(I_k)
=
\bigl((\grad P_0(-),a)\mid a\in\g\bigr)
=
J_{P_0}
=
J_P.
\]
Finally, passing to the vertex algebra quotient gives
\[
R_{L_k(\g)}
=
R_{V^k(\g)}/\sigma(I_k).
\]
All maps preserve both the grading and the Poisson bracket, which
proves the claim.
\end{proof}

\begin{remark}
The boundary admissible condition is used precisely to obtain
$E\cong\g$. At a general simply-laced principal admissible level
$k+h^\vee=a/b$, with $\gcd(a,b)=1$ and $a\geq h^\vee$, the
corresponding generating space has highest weight
$(a-h^\vee+1)\theta$. It need not be the adjoint representation, so
the symmetry used in \eqref{eq:closed} does not follow from the same
argument.
\end{remark}

\section{Restriction to Slodowy slices}\label{sec:slices}
The following restriction property is established in
\cite[Theorem~1.4]{Kang}. We recall a proof for completeness.

For a nonzero nilpotent element $f$, choose an
$\mathfrak{sl}_2$-triple $(e,h,f)$ and set
\[
\mathcal S_f=f+\g^e.
\]
For $f=0$, set $\mathcal S_0=\g$.

\begin{lemma}\label{lem:slice}
For every $G$-invariant polynomial $P\in\C[\g]^G$, one has
\[
J_P\C[\mathcal S_f]
=
J_{P|_{\mathcal S_f}}.
\]
\end{lemma}

\begin{proof}
The assertion is immediate for $f=0$, so assume that $f\neq0$.
Put $A=\C[\mathcal S_f]$. The transversality of the Slodowy slice
gives
\[
[\g,x]+\g^e=\g
\qquad
\text{for every }x\in\mathcal S_f.
\]
Consider the $A$-linear map
\[
\Phi:A\otimes(\g\oplus\g^e)\longrightarrow A\otimes\g,
\qquad
(a(x),b(x))\longmapsto[a(x),x]+b(x).
\]
Its cokernel has zero fiber at every closed point of
$\mathcal S_f$. Since the cokernel is a finitely generated
$A$-module, its localization at every maximal ideal vanishes by
Nakayama's lemma. Hence $\Phi$ is surjective.

For any constant $u\in\g$, choose polynomial maps $a$ and $b$ such
that
\[
u=[a(x),x]+b(x).
\]
The invariance of $P$ gives
\[
(\partial_uP)(x)
=
dP_x(b(x))
=
d(P|_{\mathcal S_f})_x(b(x)).
\]
Expressing $b$ in a fixed basis of $\g^e$ proves
\[
J_PA\subset J_{P|_{\mathcal S_f}}.
\]
Conversely, every constant tangent vector to
$\mathcal S_f=f+\g^e$ is given by an element of $\g^e\subset\g$.
Hence every directional derivative of $P|_{\mathcal S_f}$ is the
restriction to $\mathcal S_f$ of the corresponding directional
derivative of $P$. Thus
\[
J_{P|_{\mathcal S_f}}\subset J_PA,
\]
which proves the assertion.
\end{proof}

\begin{proof}[Proof of \Cref{th:main2}]
By \eqref{eq:DS-scheme} and \Cref{th:main1},
\[
R_{H_{DS,f}^0(L_k(\g))}
\cong
\C[\mathcal S_f]/J_P\C[\mathcal S_f].
\]
Applying \Cref{lem:slice} gives the asserted algebra isomorphism. The
compatibility of the Poisson structures follows from
Drinfeld--Sokolov reduction; in particular, the resulting Jacobian
ideal is a Poisson ideal for the Slodowy bracket.
\end{proof}

\section{Non-simply-laced principal boundary levels}\label{sec:nonsimply}
We now allow $\g$ to be non-simply-laced. Normalize the invariant
bilinear form so that long roots have squared length two, and let
$r^\vee$ denote the lacing number. For $\alpha\in\h^*$, let
$\alpha^\sharp\in\h$ denote the image of $\alpha$ under the
identification $\h^*\cong\h$ induced by the invariant form. Thus
\[
\alpha^\vee=\frac{2\alpha^\sharp}{(\alpha,\alpha)}.
\]
Assume that
\begin{equation}\label{eq:non-level}
k+h^\vee=\frac{h^\vee}{p},\qquad
p\geq2,\qquad
\gcd(p,h^\vee)=\gcd(p,r^\vee)=1.
\end{equation}

\begin{theorem}\label{th:non}
Let $\g$ be a non-simply-laced simple Lie algebra, and let $k$ and
$p$ satisfy \eqref{eq:non-level}. Then there exists a nonzero
homogeneous invariant polynomial $P\in\C[\g]^G$ of degree $p+1$ such
that
\[
R_{L_k(\g)}\cong\C[\g]/J_P
\]
as graded Poisson algebras.
Moreover, its restriction $q=P|_{\h}$ satisfies
\[
D_u(\partial_v q)=0
\qquad (u,v\in\h),
\]
where
\[
D_u
=
\partial_u-
\sum_{\alpha>0}
c_\alpha
\frac{\alpha(u)}{\alpha(x)}(1-s_\alpha),
\qquad
c_\alpha=\frac{p(\alpha,\alpha)}{2h^\vee}.
\]
Equivalently,
\[
c_{\mathrm{long}}=\frac{p}{h^\vee},
\qquad
c_{\mathrm{short}}=\frac{p}{r^\vee h^\vee}.
\]
Here $s_\alpha$ acts on scalar polynomials by
$(s_\alpha\phi)(x)=\phi(s_\alpha x)$. Thus the displayed condition
means that every component of $\grad q$ is a singular polynomial.
The conclusion of \Cref{th:main2} also holds for this potential $P$.
\end{theorem}
\begin{proof}
Put $\lambda=k\Lambda_0$. By \cite[Proposition~2.1]{KW3}, the
additional integral simple root is $\beta=p\delta-\theta$. Since
$\theta$ is long, $\beta^\vee=pK-\theta^\vee$, and hence
\[
\langle\lambda+\widehat\rho,\beta^\vee\rangle
=p(k+h^\vee)-(h^\vee-1)=1.
\]
Thus Steps~1 and~2 of the proof of \Cref{th:main1} apply unchanged:
the generating singular space is an adjoint module in degree $p$,
and its nonzero symbol defines a $G$-equivariant polynomial map
$F:\g\to\g$ whose components generate $\sigma(I_k)$.

Choose root vectors with $(e_\alpha,f_\alpha)=1$ and
$[e_\alpha,f_\alpha]=\alpha^\sharp$. The radial identity gives
\[
\partial_uF(x)=\frac{p}{h^\vee}\sum_{\alpha>0}
\frac{\alpha(u)\alpha(F(x))}{\alpha(x)}\alpha^\sharp,
\qquad x\in\h_{\mathrm{reg}},\quad u\in\h.
\]
Pairing with $v\in\h$ makes the right-hand side symmetric in $u,v$.
As in Step~3, the polynomial
$q(x)=(x,F(x))/(p+1)$ is therefore nonzero and $W$-invariant, with
$\grad q=F|_\h$. Its Chevalley lift $P$ satisfies $\grad P=F$ by
the density argument of Step~4. Consequently, $\sigma(I_k)=J_P$.

For $\phi_v(x)=(F(x),v)=\partial_vq(x)$, equivariance gives
\[
\phi_v(x)-\phi_v(s_\alpha x)
=\frac{2\alpha(v)\alpha(F(x))}{(\alpha,\alpha)}.
\]
Together with the radial identity, this proves $D_u\phi_v=0$ for the
stated multiplicities. Finally, the assertion about
Drinfeld--Sokolov reduction follows from \Cref{lem:slice} and
\eqref{eq:DS-scheme}.
\end{proof}

\begin{remark}
In types $B_n$ and $C_n$, explicit singular polynomials spanning
a copy of the reflection representation are given by
Chmutova--Etingof \cite[Proposition~4.1]{CE}, as a special case
of their construction for $G(\ell,1,n)$ with $\ell=2$.
\end{remark}

\end{document}